\documentclass[a4paper,12pt]{amsart}
\usepackage[margin=3cm]{geometry}
\newtheorem{theorem}{Theorem}[section]
\newtheorem{rem} [theorem] {Remark}
\newtheorem{prop} [theorem]{Proposition}
\newtheorem{lemma}[theorem]{Lemma}
\newtheorem{definition}[theorem]{Definition}
\newcommand{\ovprt}{\overline{\partial}}

\newcommand{\dquer}{\overline\partial}
\newcommand{\dquers}{\overline\partial ^*_\varphi}
\newcommand{\boxphi}{\square_\varphi}
\newcommand{\levim}{\frac{\partial^2\varphi}{\partial z_j\partial\overline z_k}}

\numberwithin{equation}{section}
\begin{document}
\title{Compactness  for  the  $\ovprt $ - Neumann problem  - a functional analysis approach.}

\author{ Friedrich Haslinger}

\thanks{Partially supported by the FWF-grant  P19147.}

 \address{ F. Haslinger: Institut f\"ur Mathematik, Universit\"at Wien,
Nordbergstrasse 15, A-1090 Wien, Austria}
\email{ friedrich.haslinger@univie.ac.at}

\maketitle

\begin{abstract} ~\\
 We discuss compactness of the $\ovprt $-Neumann operator in the setting of weighted $L^2$-spaces on $\mathbb{C}^n.$
For this purpose we use  a description of relatively compact subsets of $L^2$- spaces.
We also point out how to use this method to show that property (P) implies compactness for the  $\ovprt $-Neumann operator on a smoothly bounded pseudoconvex domain and mention an abstract functional analysis characterization of compactness of the $\ovprt $-Neumann operator. 

\end{abstract}

\section{Introduction.}~\\

 In this paper we continue the investigations of \cite{HaHe} concerning  existence and compactness of the canonical solution operator to $\ovprt $ on weighted
$L^2$-spaces over $\mathbb C^n.$

Let $\varphi : \mathbb C^n \longrightarrow \mathbb R^+ $ be a plurisubharmonic $\mathcal C^2$-weight function and define the space
$$L^2(\mathbb C^n , \varphi )=\{ f:\mathbb C^n \longrightarrow \mathbb C \ : \ \int_{\mathbb C^n}
|f|^2\, e^{-\varphi}\,d\lambda < \infty \},$$
where $\lambda$ denotes the Lebesgue measure, the space $L^2_{(0,1)}(\mathbb C^n, \varphi )$ of $(0,1)$-forms with coefficients in
$L^2(\mathbb C^n , \varphi )$ and the space $L^2_{(0,2)}(\mathbb C^n, \varphi )$ of $(0,2)$-forms with coefficients in
$L^2(\mathbb C^n , \varphi ).$
Let 
$$\langle f,g\rangle_\varphi=\int_{\mathbb{C}^n}f \,\overline{g} e^{-\varphi}\,d\lambda$$
denote the inner product and 
$$\| f\|^2_\varphi =\int_{\mathbb{C}^n}|f|^2e^{-\varphi}\,d\lambda $$
the norm in $L^2(\mathbb C^n , \varphi ).$

We consider the weighted
$\ovprt $-complex 
$$
L^2(\mathbb C^n , \varphi )\underset{\underset{\ovprt_\varphi^* }
\longleftarrow}{\overset{\ovprt }
{\longrightarrow}} L^2_{(0,1)}(\mathbb C^n , \varphi )
\underset{\underset{\ovprt_\varphi^* }
\longleftarrow}{\overset{\ovprt }
{\longrightarrow}} L^2_{(0,2)}(\mathbb C^n , \varphi ),
$$
where $\ovprt_\varphi^*$ is the adjoint operator to $\ovprt $ with respect to the weighted inner product. For $u=\sum_{j=1}^nu_jd\overline z_j\in {\text {dom}}(\ovprt_\varphi^*)$ one has
$$\ovprt_\varphi^*u=-\sum_{j=1}^n \left ( \frac{\partial}{\partial z_j}-
\frac{\partial \varphi}{\partial z_j}\right )u_{j}.$$
The complex Laplacian on $(0,1)$-forms is defined as
$$\boxphi := \dquer  \,\dquers + \dquers \dquer,$$
where the symbol $\boxphi $ is to be understood as the maximal closure of the operator initially defined on forms with coefficients in $\mathcal{C}_0^\infty$, i.e., the space of smooth functions with compact support.

$\boxphi $ is a selfadjoint and positive operator, which means that 
$$\langle \boxphi f,f\rangle_\varphi \ge 0 \ , \   {\text{for}} \  f\in dom (\boxphi ).$$
The associated Dirichlet form is denoted by 
\begin{equation}\label{diri}
Q_\varphi (f,g)= \langle \dquer f,\dquer g\rangle_\varphi + \langle \dquers f ,\dquers g\rangle_\varphi, 
\end{equation}
for $f,g\in dom (\dquer ) \cap dom (\dquers ).$ The weighted $\dquer $-Neumann operator
$N_\varphi $ is - if it exists - the bounded inverse of $\boxphi .$ 
\vskip 0.5 cm
We indicate that  $f\in dom(\dquers)$ if and only if 
\begin{equation*}
 \sum_{j=1}^n\left ( \frac{\partial f_j}{\partial z_j}- \frac{\partial \varphi}{\partial z_j}\, f_j \right )
 \in L^2(\mathbb{C}^n, \varphi )
 \end{equation*}
  and that forms with coefficients in $\mathcal{C}_0^\infty(\mathbb{C}^n)$ are dense in $dom(\dquer)\cap dom(\dquers)$ in the graph norm $f\mapsto (\Vert \dquer f\Vert _\varphi^2+\Vert \dquers f\Vert _\varphi^2)^\frac{1}{2}$ (see \cite{GaHa}).
\vskip 0.5 cm
Now we suppose that the lowest eigenvalue $\mu_\varphi $ of the Levi - matrix 
$$M_\varphi=\left(\levim\right)_{jk}$$
of $\varphi $ satisfies
$$\liminf_{|z|\to \infty}\mu_\varphi (z)>0, \ \ \ (^*)$$
and mention the Kohn-Morrey formula:
\begin{equation}\label{kohn}
\| \ovprt u \|^2_\varphi + \| \ovprt_\varphi ^* u\|^2_\varphi 
=\sum_{j,k =1}^n \int_{\mathbb{C}^n} \left |\frac{\partial u_j}{\partial \overline{z}_k}
\right |^2 \,e^{-\varphi}\,d\lambda + \int_{\mathbb{C}^n} \sum_{j,k=1}^n\frac{\partial^2 \varphi}{\partial z_j \partial\overline{z}_k}\,u_j\overline{u}_k\, e^{-\varphi}\,d\lambda
\end{equation}
from which we get
\begin{equation}\label{komo}
\int_{\mathbb{C}^n} \sum_{j,k=1}^n\frac{\partial^2 \varphi}{\partial z_j \partial\overline{z}_k}\,u_j\overline{u}_k\, e^{-\varphi}\,d\lambda \le \|\overline{\partial}u\|^2_\varphi +
\| \ovprt_\varphi ^* u\|^2_\varphi ,
\end{equation}

hence for a  plurisubharmonic weight function $\varphi $ satisfying (*), there is a $C>0$ such that
$$\|u\|^2_\varphi \le C( \| \ovprt u \|^2_\varphi + \| \ovprt_\varphi ^* u\|^2_\varphi )$$
for each $(0,1)$-form $u\in $dom\,$(\ovprt ) \,\cap$
dom\,$(\ovprt_\varphi^*).$

For the proof see \cite{FS}, \cite{GaHa} or \cite{Str}.

Now it follows that there exists a uniquely determined $(0,1)$-form 
\newline $N_\varphi u\in $dom\,($\ovprt) \,\cap$ dom\,($\ovprt_\varphi^*)$ such that
$$\langle u,v\rangle_\varphi=Q_\varphi (N_\varphi u,v)= \langle\ovprt N_\varphi u,\ovprt v\rangle_\varphi + \langle\ovprt_\varphi^* N_\varphi u,
\ovprt_\varphi^* v\rangle_\varphi,$$
and  that
\begin{equation}\label{graph}
\|\ovprt N_\varphi u \|_\varphi ^2 + \|\ovprt_\varphi^* N_\varphi u \|_\varphi ^2 \le  C_1 \|u\|_\varphi ^2
\end{equation}
which means that 
$$N_{1,\varphi}: L^2_{(0,1)}(\mathbb{C}^n, \varphi ) \longrightarrow {\text {dom}}\,(\ovprt )\,\cap
{\text {dom}}\,(\ovprt_\varphi^*)$$
is continuous in the graph topology,
as well as
$$\|N_\varphi u\|_\varphi ^2 \le C_2 (\|\ovprt N_\varphi u \|_\varphi ^2 + \|\ovprt_\varphi^* N_\varphi u \|_\varphi ^2) \le  C_3 \|u\|_\varphi ^2, $$
where $C_1, C_2, C_3>0$ are constants. Hence we get that $N_\varphi$ is a continuous linear operator from $L^2_{(0,1)}( \mathbb{C}^n, \varphi)$ into itself (see also \cite{ChSh}). 

We will give a new proof of the main result in \cite{HaHe} using a direct approach, see \cite{B}, Corollaire IV.26, where two conditions are given which imply that a subset of an $L^2$- space is relatively compact. The first of these conditions will correspond to G\aa rding's inequality (see for instance \cite{F} , \cite{GaHa},) and the second condition corresponds to our assumption on the lowest eigenvalue of the Levi matrix $M_\varphi .$

 We indicate  how to use this method to show that property (P)  implies compactness for the  $\ovprt $-Neumann operator on a smoothly bounded pseudoconvex domain $\Omega\subset\subset \mathbb{C}^n$ and finally mention an abstract necessary and sufficient condition for the $\ovprt $- Neumann operator to be compact.

\vskip 1 cm

\section{Weighted Sobolev spaces}~\\

Now we define an appropriate Sobolev space and prove compactness of the corresponding embedding, for related settings see \cite{BDH}, \cite{Jo}, \cite{KM} .

\begin{definition} Let
$$\mathcal{W}^{Q_\varphi}= \{ u\in L^2_{(0,1)}(\mathbb{C}^n , \varphi) \ : \  \| \ovprt u \|^2_\varphi + \| \ovprt_\varphi ^* u\|^2_\varphi< \infty \}$$ with norm
$$\| u\|_{Q_\varphi}  =  (\| \ovprt u \|^2_\varphi + \| \ovprt_\varphi ^* u\|^2_\varphi )^{1/2}.$$
\end{definition}

{\bf Remark:} $\mathcal{W}^{Q_\varphi}$ coincides with the form domain $dom (\dquer ) \cap dom (\dquers )$ of $Q_\varphi $ (see \cite{Ga}, \cite{GaHa} ).

\begin{prop}\label{embed}Suppose that the weight function  $\varphi$ is  plurisubharmonic  and that the lowest eigenvalue $\mu_{\varphi}$ of the Levi - matrix $M_{\varphi }$ satisfies
$$\lim_{|z|\rightarrow \infty}\mu_\varphi(z)  = +\infty\, . \ \ \ (^{**})$$
Then the embedding 
$$j_\varphi : \mathcal{W}^{Q_\varphi } \hookrightarrow  L^2_{(0,1)}(\mathbb{C}^n, \varphi )$$ is compact.
\end{prop}
\begin{proof}
For $u\in \mathcal{W}^{Q_\varphi}$ we have by \ref{komo}
$$ \| \ovprt u \|^2_\varphi + \| \ovprt_\varphi ^* u\|^2_\varphi \ge  \langle M_{\varphi } u, u\rangle_\varphi .$$

This implies
\begin{equation}\label{compact}
 \| \ovprt u \|^2_\varphi + \| \ovprt_\varphi ^* u\|^2_\varphi \ge  \int_{\mathbb{C}^n}\mu_\varphi(z) \, |u(z)|^2\,  e^{-\varphi(z)}\,d\lambda (z) .
\end{equation}
\vskip 0.5 cm

We show that the unit ball in $\mathcal{W}^{Q_\varphi }$ is relatively compact in $ L^2_{(0,1)}(\mathbb{C}^n, \varphi ).$
For this purpose we use the following lemma, see for instance \cite{B} Corollaire IV.26.
\vskip 0.5 cm
\begin{lemma}\label{brezis} Let $\mathcal{A}$ be a bounded subset 
of $ L^2(\mathbb{C}^n, \varphi ).$ Suppose that

(i) for each $\epsilon >0$ and for each $R>0$ there exists $\delta >0$ such that 
$$ \| \tau_h f -f \|_{ L^2(\mathbb{B}_R, \varphi )} < \epsilon $$
for each $h\in \mathbb{C}^n$ with $|h|<\delta $ and for each $f\in \mathcal{A},$ where $\tau_hf(z)=f(z+h)$ and $\mathbb{B}_R=\{ z \in \mathbb{C}^n : |z|<R \};$

(ii) for each $\epsilon >0$ there exists $R>0$ such that 

$$\|f\|_{ L^2(\mathbb{C}^n \setminus \mathbb{B}_R, \varphi )} < \epsilon $$
for each $f\in \mathcal{A}.$

Then $\mathcal{A}$ is relatively compact in $ L^2(\mathbb{C}^n, \varphi ).$
\end{lemma}

\begin{rem} Conditions (i) and (ii) are also necessary for $\mathcal{A}$ to be relatively compact in $ L^2(\mathbb{C}^n, \varphi )$ \ (see \cite{B}).
\end{rem}

First we show that condition (i) of Lemma \ref{brezis} is satisfied in our situation. 
Let $u=\sum_{j=1}^n u_j\, dz_j$ be a $(0,1)$-form with coefficients in $\mathcal{C}_0^\infty.$
For each $u_j$ and for $t\in \mathbb{R}$ and $h=(h_1,\dots , h_n)\in \mathbb{C}^n$ let 
$$v_j(t):= u_j(z+th).$$
Note that 
$$|v_j'(t)|\le |h| \left [ \sum_{k=1}^n\left ( \left | \frac{\partial u_j}{\partial x_k}(z+th)\right |^2+
\left | \frac{\partial u_j}{\partial y_k}(z+th)\right |^2\right ) \right ]^{1/2},$$
where $z_k=x_k+iy_k,$ for $k=1,\dots ,n.$
By the fact that 
$$u_j(z+h)-u_j(z)=v_j(1)-v_j(0)=\int_0^1 v_j'(t)\,dt  $$
 we can now estimate for $|h|<R$
$$\int_{\mathbb{B}_R}|\tau_hu_j(z)-u_j(z)|^2 e^{-\varphi (z)}\,d\lambda (z) 
=\int_{\mathbb{B}_R}|\tau_h(\chi_R u_j)(z)-\chi_Ru_j(z)|^2 e^{-\varphi (z)}\,d\lambda (z) $$
$$\le |h|^2 
\,  \int_{\mathbb{B}_{R}}\left [ \int_0^1 \sum_{k=1}^n\left ( \left | \frac{\partial (\chi_Ru_j)}{\partial x_k}(z+th)\right |^2+
\left | \frac{\partial (\chi_Ru_j)}{\partial y_k}(z+th)\right |^2\right )\,dt \, \right ]
 e^{-\varphi (z)} \,d\lambda (z) 
$$ 
$$ \le C_{R,\varphi}\,  |h|^2 
\,  \int_{\mathbb{B}_{3R}} \sum_{k=1}^n\left ( \left | \frac{\partial (\chi_Ru_j)}{\partial x_k}(z)\right |^2+
\left | \frac{\partial (\chi_Ru_j)}{\partial y_k}(z)\right |^2\right )
 e^{-\varphi (z)} \,d\lambda (z) 
$$
\vskip 0.3 cm
for $j=1,\dots ,n$ where $\chi_R$ is a $\mathcal{C}^\infty $ cutoff function which is identically $1$ on $\mathbb{B}_{2R}$ and zero outside  $\mathbb{B}_{3R}$
and by 
 G\aa rding's inequality for $\mathbb{B}_{3R}$ (see \cite{ChSh}, \cite{F}, \cite{GaHa})
$$
\begin{array}{ll}
\Vert \chi_R u \Vert ^2_{\varphi, 1}&\leq C'_{\varphi,R}
\left(\Vert\dquer ( \chi_R u)\Vert^2_\varphi+\Vert\dquer^*_\varphi (\chi_R u)\Vert^2_\varphi +\Vert \chi_R u\Vert^2_\varphi \right) \\
&\leq C''_{\varphi,R} \left(\Vert\dquer  u\Vert^2_\varphi +\Vert\dquer^*_\varphi  u\Vert^2_\varphi +\Vert  u\Vert^2_\varphi \right)
\end{array}$$
we can control the last integral by the norm $\|u\|^2_{Q_\varphi}.$
Since we started from the unit ball in $\mathcal{W}^{Q_\varphi }$ we get that condition (i) of Lemma \ref{brezis} is satisfied.

\vskip 0.3 cm
Condition (ii) of Lemma \ref{brezis} is satisfied for the unit ball of $\mathcal{W}^{Q_\varphi }$ since we have
\vskip 0.3 cm

$$\int_{\mathbb{C}^n \setminus \mathbb{B}_R} |u(z)|^2 e^{-\varphi (z)}\,d\lambda (z)\le
\int_{\mathbb{C}^n \setminus \mathbb{B}_R} \frac{\mu_\varphi(z) |u_(z)|^2}{\inf \{\mu_\varphi(z)  :  |z|\ge R\}}   e^{-\varphi(z)}  d\lambda (z).$$
\vskip 0.1 cm

So formula (\ref{compact}) together with assumption (**) shows that

$$\int_{\mathbb{C}^n \setminus \mathbb{B}_R} |u(z)|^2 e^{-\varphi (z)}\,d\lambda (z)\le
 \frac{ \|u\|^2_{Q_\varphi}}{\inf \{\mu_\varphi(z) \, : \, |z|\ge R\}} < \epsilon,$$

 if $R$ is big enough.

\end{proof}

We are now able to give a short proof of the main result in \cite{HaHe} or \cite{GaHa}
\begin{prop}
\label{cp 1}
Let $\varphi$ be a plurisubharmonic $\mathcal C^2$- weight function. If the lowest eigenvalue $\mu_\varphi (z)$ of the Levi - matrix $M_\varphi$ satisfies $(^{**})$, then $N_\varphi$ is compact.
\end{prop}

\begin{proof}  By Proposition \ref{embed},   the embedding $\mathcal{W}^{Q_\varphi} \hookrightarrow  L^2_{(0,1)}(\mathbb{C}^n , \varphi)$ is compact. The inverse $N_\varphi $ of  $\square_{\varphi }$ is 
continuous as an operator from $L^2_{(0,1)}(\mathbb{C}^n , \varphi)$ into $\mathcal{W}^{Q_\varphi},$
this follows from \ref{graph}. Therefore we have that $N_\varphi $ is compact as an operator from $ L^2_{(0,1)}(\mathbb{C}^n , \varphi)$ into itself.

\end{proof}

Now notice that 
$$N_\varphi : L^2_{(0,1)}(\mathbb{C}^n , \varphi) \longrightarrow L^2_{(0,1)}(\mathbb{C}^n , \varphi) $$
can be written in the form 
$$N_\varphi = j_\varphi \circ j_\varphi ^* \ , $$
where 
$$j_\varphi ^* : L^2_{(0,1)}(\mathbb{C}^n , \varphi) \longrightarrow \mathcal{W}^{Q_\varphi}$$
is the adjoint operator to $j_\varphi$  \ (see \cite{Str}).

This means that $N_\varphi $ is compact if and only if $j_\varphi$ is compact and summarizing the above results we get the following

\begin{prop} Let $\varphi : \mathbb C^n \longrightarrow \mathbb R^+ $ be a plurisubharmonic $\mathcal C^2$-weight function . The $\ovprt $-Neumann operator 
$$N_\varphi  :  L^2_{(0,1)}(\mathbb{C}^n , \varphi) \longrightarrow L^2_{(0,1)}(\mathbb{C}^n , \varphi) $$ 
is compact if and only if  for each $\epsilon >0$ there exists $R>0$ such that 

$$\|u\|_{ L^2_{(0,1)}(\mathbb{C}^n \setminus \mathbb{B}_R, \varphi )} < \epsilon $$
for each $u\in \mathcal{W}^{Q_\varphi}$ with
$$\| \ovprt u \|^2_\varphi + \| \ovprt_\varphi ^* u\|^2_\varphi \le 1.$$

\end{prop}

\vskip 1 cm

\section{Smoothly bounded pseudoconvex domains and \\
properties (P) and (\~P)}~\\

\vskip 0.2 cm

Let $\Omega \subset\subset \mathbb{C}^n$ be a smoothly bounded pseudoconvex domain. $\Omega $ satisfies property (P), if or each $M>0$ there exists a a neighborhood $U$ of $\partial \Omega$ and a plurisubharmonic function $\varphi_M \in \mathcal{C}^2 (U)$ such that 
$$\sum_{j,k=1}^n \frac{\partial^2 \varphi_M}{\partial z_j \partial \overline z_k}(p) t_j \overline t_k \ge M \| t \|^2,$$
for all $p\in \partial \Omega$ and for all $t\in \mathbb{C}^n.$
\vskip 0.2 cm
$\Omega $ satisfies property (\~P) if the following holds: there is a constant $C$ such that for all $M>0$ there exists a $\mathcal{C}^2$ function $\varphi_M$ in a neighborhood $U$ (depending on $M$) of $\partial \Omega$ with

(i) $  \left | \sum_{j=1}^n \frac{\partial \varphi_M}{\partial z_j}(z) t_j \right |^2 \le C \sum_{j=1}^n \frac{\partial^2 \varphi_M}{\partial z_j \partial \overline z_k}(z) t_j \overline t_k$

and

(ii) $\sum_{j=1}^n \frac{\partial^2 \varphi_M}{\partial z_j \partial \overline z_k}(z) t_j \overline t_k \ge M \| t \|^2,$
 
for all $z\in U$ and for all $t\in \mathbb{C}^n.$
\vskip 0.2 cm

In \cite{C} Catlin showed that condition (P) implies compactness of the $\overline \partial$- operator $N$ on $L^2_{(0,1)}(\Omega )$ and McNeal (\cite{McN}) showed that property (\~P) also implies  compactness of the $\overline \partial$- operator $N$ on $L^2_{(0,1)}(\Omega ).$ It is not difficult to show that property (P) implies property (\~P), see for instance \cite{Str}.

We can now use a similar approach as in section 2 to prove Catlin's result. For this purpose we use the following version of  lemma \ref{brezis}

\begin{lemma}\label{brezis1} Let $\mathcal{A}$ be a bounded subset 
of $ L^2(\Omega).$ Suppose that

(i) for each $\epsilon >0$ and for each $\omega \subset\subset \Omega$ there exists $\delta >0, \delta < {\text{dist}}(\omega, \Omega^c)$ such that 
$$ \| \tau_h f -f \|_{ L^2(\omega )} < \epsilon $$
for each $h\in \mathbb{C}^n$ with $|h|<\delta $ and for each $f\in \mathcal{A}$,

(ii) for each $\epsilon >0$ there exists $\omega \subset\subset  \Omega$ such that 

$$\|f\|_{ L^2(\Omega \setminus \omega)} < \epsilon $$
for each $f\in \mathcal{A}.$

Then $\mathcal{A}$ is relatively compact in $ L^2(\Omega ).$
\end{lemma}

\begin{rem} Conditions (i) and (ii) are also necessary for $\mathcal{A}$ to be relatively compact in $ L^2(\Omega ).$
\end{rem}

In order to show that the unit ball in 
$dom(\dquer)\cap dom(\dquer^*)$ in the graph norm $f\mapsto (\Vert \dquer f\Vert ^2+\Vert \dquer^* f\Vert^2)^\frac{1}{2}$
 satisfies condition (i) of \ref{brezis1} we remark that
 G\aa rding's inequality holds for $\omega \subset\subset \Omega$ (see section 2). To verify condition (ii) we use property (P) and the following version of the Kohn-Morrey formula
 
\begin{equation}\label{komo1}
\int_{\Omega} \sum_{j,k=1}^n\frac{\partial^2 \varphi_M}{\partial z_j \partial\overline{z}_k}\,u_j\overline{u}_k\, e^{-\varphi_M}\,d\lambda \le \|\overline{\partial}u\|^2_{\varphi_M}+
\| \ovprt_{\varphi_M} ^* u\|^2_{\varphi_M} ,
\end{equation}
here we used  that $\Omega $ is pseudoconvex, which means that the boundary terms in the Kohn-Morrey formula can be neglected. Now we point out that the weighted $\overline \partial $- complex is equivalent to the unweighted one and that the expression $\sum_{j=1}^n \frac{\partial \varphi_M}{\partial z_j}u_j$ which appears in $\ovprt_{\varphi_M} ^* u,$ can be controlled by the complex Hessian  $\sum_{j,k=1}^n\frac{\partial^2 \varphi_M}{\partial z_j \partial\overline{z}_k}\,u_j\overline{u}_k ,$ which follows from the fact that property (P) implies property (\~{P}) (see \cite{Str}). Of course we also use that the weight $\varphi_M$ is bounded on $\Omega \subset\subset \mathbb{C}^n.$ In this way the same reasoning as in section 2 shows that property (P) implies condition (ii) of lemma \ref{brezis1}. Therefore condition
(P) gives that the unit ball of
$dom(\dquer)\cap dom(\dquer^*)$ in the graph norm $f\mapsto (\Vert \dquer f\Vert ^2+\Vert \dquer^* f\Vert^2)^\frac{1}{2}$
 is relatively compact in $L^2_{(0,1)}(\Omega)$ and hence that the $\overline \partial $-Neumann operator is compact.

\vskip 0.3 cm
Now let 
$$j: dom(\dquer)\cap dom(\dquer^*) \hookrightarrow L^2_{(0,1)}(\Omega)$$
denote the  embedding. It follows from \cite{Str} that 
$$N= j \circ j^*.$$
Hence $N$ is compact if and only if $j$ is compact, where $dom(\dquer)\cap dom(\dquer^*)$ is endowed with the graph norm  $f\mapsto (\Vert \dquer f\Vert ^2+\Vert \dquer^* f\Vert^2)^\frac{1}{2}.$
 
 \begin{prop} Let $\Omega \subset\subset \mathbb{C}^n$ be a smoothly bounded pseudoconvex domain. Let $\mathcal{B}$ denote the unit ball of $dom(\dquer)\cap dom(\dquer^*)$ in the graph norm $f\mapsto (\Vert \dquer f\Vert ^2+\Vert \dquer^* f\Vert^2)^\frac{1}{2}.$
 
 The $\ovprt $- Neumann operator $N$ is compact if and only if $\mathcal{B}$ as a subset of $L^2_{(0,1)}(\Omega)$ satisfies the following condition:
 
   for each $\epsilon >0$ there exists $\omega \subset\subset  \Omega$ such that 

$$\|f\|_{ L^2_{(0,1)}(\Omega \setminus \omega)} < \epsilon $$
for each $f\in \mathcal{B}.$
\end{prop}

This follows from the above remarks about the embedding $j$ and the fact that the two conditions in \ref{brezis1} are also necessary for a bounded set in $L^2$ to be relatively compact. For a localized version of the above result see \cite{Sa}.
\vskip 0.5 cm
\textbf{Acknowledgement. } The  author is grateful to the Erwin Schr\"odinger Institute (ESI), Vienna, where most of this work was carried out and to Mehmet Celik, Bernard Helffer and Emil Straube for valuable discussions during the ESI - program on the $\ovprt $-Neumann operator.


\vskip 1 cm
 \begin{thebibliography}{ABCD}


\vskip 0.5 cm

\bibitem[BDH]{BDH} P. Bolley, M. Dauge and B. Helffer,{\em Conditions suffisantes pour l'injection compacte d'espace de Sobolev \`a poids,} S\'eminaire \'equation aux d\'eriv\'ees partielles (France), vol.1, Universit\'e de Nantes (1989), 1--14.

\bibitem[B]{B} H. Brezis, {\em Analyse fonctionnelle, Th\'eorie et applications}, Masson , Paris, 1983.

\bibitem[C]{C} D.W. Catlin, {\em Global regularity of the $\overline \partial $-Neumann operator,} Proc. Symp. Pure Math. {\bf 41} (1984), 39-49.
\bibitem[ChSh]{ChSh} So-Chin Chen and Mei-Chi Shaw, {\em Partial differential equations in several complex variables}, Studies in Advanced Mathematics, Vol.~19, Amer. Math. Soc. 2001.
\bibitem[F]{F} G.B. Folland, {\em Introduction to partial differential equations}, Princeton University Press, Princeton, 1995.

\bibitem[FS]{FS}  S.~Fu and E.J.~Straube, {\em Compactness in
the $\ovprt -$Neumann problem}, Complex Analysis and Geometry (J. McNeal, ed.),
Ohio State Math. Res. Inst. Publ. {\bf 9} (2001), 141--160.



\bibitem[Ga]{Ga} K. Gansberger, Compactness of the $\dquer$-Neumann operator, Dissertation, University of Vienna, 2009.
\bibitem[GaHa]{GaHa} K. Gansberger and F. Haslinger, {\em Compactness estimates for the $\overline \partial$- Neumann problem in weighted $L^2$- spaces}, Proceedings of the conference on Complex Analysis 2008 in honour of Linda Rothschild, Fribourg 2008, to appear.
\bibitem[HaHe]{HaHe} F. Haslinger and B. Helffer, {\em Compactness of the solution operator to $\ovprt $ in weighted $L^ 2$ - spaces}, J. of Functional Analysis,  {\bf 243} (2007), 679-697.
 



\bibitem[Jo] {Jo} J.~Johnsen,
{\em  On the spectral properties of Witten Laplacians, their range
  projections
 and Brascamp-Lieb's inequality}, Integral Equations Operator Theory
{\bf 36}  (3), 2000, 288--324.

\bibitem[KM]{KM} J.-M. Kneib and F. Mignot, {\em Equation de Schmoluchowski g\'en\'eralis\'ee,} Ann. Math. Pura Appl. (IV) {\bf 167} (1994), 257--298.


\bibitem[McN]{McN} J.D. McNeal , {\em A sufficient condition for compactness of the $\ovprt $-Neumann operator}, J. of Functional Analysis,  {\bf 195} (2002), 190-205.
\bibitem[Sa]{Sa} S. Sahutoglu , {\em Compactness of the $\ovprt $- Neumann problem and Stein neighborhood bases}, Dissertation Texas A \& M University, 2006.

\bibitem[Str]{Str} E. Straube, {\em The $L^2$-Sobolev theory of the $\ovprt $-Neumann problem,} ESI Lectures in Mathematics and Physics, EMS (to appear).

\end{thebibliography}
\end{document}